\magnification = 1200
\showboxdepth=0 \showboxbreadth=0

\baselineskip 14pt
\parskip3pt
\def\qed{\hfill\vrule height6pt width6pt depth0pt}

\def\ss{\smallskip}
\def\ms{\medskip}
\def\bs{\bigskip}

\def\cl{\centerline}

\def\nind{\noindent}

\def\ref#1#2{\nind\hangindent.5in\hbox to .5in{#1\hfill}#2}
\def\reff#1#2{\nind\hangindent.8in\hbox to .8in{\bf #1\hfill}#2\par}
\def\refd#1#2{\nind\hangindent.8in\hbox to .8in{\bf #1\hfill {\rm--}}#2\par}
\def\pmb#1{\setbox0=\hbox{#1}
\kern-0.025em\copy0\kern-\wd0
\kern.05em\copy0\kern-\wd0
\kern-.025em\raise.0433em\box0}
\def\ca#1{{\cal #1}}

\def\Cal{\ca}

\def\frac#1#2{{#1\over#2}}

\def\text#1{\rm#1}

\outer\def\stmnt  #1. #2\par{\medbreak
\noindent{\bf#1.\enspace}{\sl#2}\par
\ifdim\lastskip<\medskipamount \removelastskip\penalty55\medskip\fi}

\def\newline{\hfill\break}
\def\:{\,:\,}

\def\({\left(}                   
\def\){\right)}

\def\[{\left[}                   
\def\]{\right]}

\def\lan{\langle}
\def\ran{\rangle}

\def\ci{\subset}

\def\fy{\infty}

\def\del{\partial}

\def\lam{\lambda}

\def\om{\omega}

\def\Om{\Omega}

\cl{\bf The comparsion principle for viscosity solutions of fully}
\cl{\bf nonlinear subelliptic equations in Carnot groups}

\bs
\cl{Changyou Wang}
\ss
\cl{Department of Mathematics, \ University of Kentucky}
\cl{Lexington, KY 40506}  
\bs
\nind{\bf Abstract}. {\it For any Carnot group $\bf G$ and a bounded domain
$\bf \Om\ci\bf G$, we prove that viscosity solutions in $C(\bar\Om)$ of
the fully nonlinear subelliptic equation $F(u,\nabla_h u,\nabla_h^2 u)=0$ are unique
when $F\in C(R\times R^m\times {\Cal S}(m))$ satisfies (i) $F$ is degenerate
subelliptic and decreasing in $u$ or (ii) $F$ is uniformly subelliptic 
and nonincreasing in $u$. This extends Jensen's uniqueness theorem from the
Euclidean space to the sub-Riemannian setting of the Carnot group.} 

\bs
\nind{\S1}. Introduction
\ss
The notion of viscosity solutions of fully nonlinear 2nd order degenerate
elliptic equation:
$$F(x,u(x),\nabla u(x), \nabla^2 u(x))=0, \hbox{ in }R^n, \eqno(1.1)$$
was developed by Crandall-Lions [CL] and Evans [E1,2] in 1980's. This idea, together
with Jensen's celebrated uniqueness theorem [J1], provides a very satisfactory theory
on existence, uniqueness, and compactness theorem of weak solutions of (1.1). The theory
of viscosity solutions has been very powerful in many applications, and we refer to 
the user's guide [CIL] by Crandall-Ishii-Lions for many such applications. 

In recent years there has been an explosion of interest in the study of analysis on
sub-Riemannian, or Carnot-Carath\'edory spaces. The corresponding developments in the
theory of partial differential equations of subelliptic type have prompted people
to consider fully nonliear equations in Carnot groups. For examples, motivated by
the very important work of Jensen [J2] on absolute minimizing Lipschitz extensions
(or ALMEs, a notion first introduced by Aronsson [A]) and viscosity solutions
to the $\fy$-laplacian equation in the Euclidean space, Bieske [B], 
Bieske-Capogna [BC], and Wang [W1] have studied absolute minimizing horizontal Lipschitz
extensions and viscosity solutions to the $\fy$-sublaplacian equation on Carnot groups. 
In particular, the notion of viscosity solutions has been extended to fully nonlinear 
subelliptic equation (see [B]) and the uniqueness of viscosity solution of $\fy$-sublaplacian
eqaution on any Carnot group was established by Wang [W1]. It is well-known (cf. 
the monographs [CC] by Caffarelli-Cabr\'e and [G] by Gutierrez) that both convexity
and the Monge-Amp\'ere equation:
$$\hbox{det}(\nabla^2 u)=f, \hbox{ in } R^n \eqno(1.2)$$
have played crucial roles in the theory of fully nonlinear elliptic equation. Inspired by
this, Lu-Manfredi-Stroffolini [LMS] and Danielli-Garofalo-Nhieu [DGN] have introduced and studied 
various notions of convexity, such as v-convexity and h-convexity, on Carnot groups (see also
[BR], [W2], [JM] for some further related results). Moreover, Garofalo-Tournier [GN] and
Gutierrez-Montanari [GM] have initiated the study of Monge-Amp\'ere measures and maximum
principle of convex functions on Heisenberg groups.   

In this paper, we are interested in the comparison principle for viscosity
solutions to 2nd order subelliptic equation which is either uniformly subelliptic, nonincreasing
or degenerate subelliptic, decreasing in the sub-Riemannian setting of the Carnot group.
In this aspect, we are able to extend Jensen's uniqueness theorem from the Eucliean space
to any Carnot group.

In order to describe our result, we first recall the basic properties of Carnot groups.
A simply connected Lie group $\bf G$ is called a Carnot group
of step $r\ge 1$, if its Lie algebra ${\it g}$ admits a vector space
decomposition in $r$ layers ${\it g}=V_1+V_2+\cdots +V_r$ such that
(i) ${\it g}$ is stratified, i.e., $[V_1, V_j]=V_{j+1}, j=1,\cdots, r-1$,
and (ii) ${\it g}$ is $r$-nilpotent, i.e. $[V_j, V_r]=0, j=1, \cdots, r$.  
We call $V_1$ the {\it horizontal} layer and $V_j, j=2,\cdots, r$
the {\it vertical} layers.
We choose an inner product $\lan\cdot,\cdot\ran$ on {\it g}
such that $V_j's$ are mutually orthogonal for $1\le j\le r$.
Let $\{X_{j, 1}, \cdots, X_{j,m_j}\}$ denote
a fixed orthonormal basis of $V_j$ for $1\le j\le r$,
where $m_j=\hbox{dim}(V_j)$ is the dimension of $V_j$.
From now on, we also denote $m=\hbox{dim}(V_1)$ as
the dimension of the horizontal layer and set
$X_i=X_{1,i}$ for $1\le i\le m$. It is well-known (see [FS])
that the exponential map ${\it exp}: {\it g}\equiv R^n\to G$ is a global
diffeomorphism and yields an exponential coordinate
system on $\bf G$, with $n=\sum_{i=i}^r m_i$ the topological dimension of $\bf G$. 
More precisely, any $p\in \bf G$ has a coordinate 
$((p_1, \cdots, p_m),(p_{2,1},\cdots, p_{2,m_2}), \cdots,
(p_{r,1},\cdots, p_{r,m_r}))$ such that
$$p={\it exp}(\xi_1(p)+\cdots \xi_r(p)), \hbox{ with }
\xi_1(p)=\sum_{l=1}^m p_lX_l,\ \xi_i(p)=\sum_{j=1}^{m_i}p_{i,j}X_{i,j}, 2\le i\le r.$$
The exponential map can induce a homogeneous pseudo-norm $N_{\bf G}$ on $\bf G$ 
in the following way (see [FS]). 
$$N_{\bf G}(p):=(\sum_{i=1}^r|\xi_i(p)|^{2r!\over i})^{1\over 2r!}, \hbox{ if }
p={\it exp}(\xi_1(p)+\cdots \xi_r(p)),\eqno(1.3)$$
where $|\xi_1(p)|=(\sum_{l=1}^m p_l^2)^{1\over 2}$, and
$|\xi_i(p)|=(\sum_{j=1}^{m_i} p_{i,j}^2)^{1\over 2}(2\le i\le r$).
Moreover, $N_{\bf G}$ yields a pseudo-distance on $\bf G$ as follows.
$$d_{\bf G}(p, q):=N_{\bf G}(p^{-1}\cdot q), \ \forall p, \ q\in\bf G, \eqno(1.4)$$
where ${\cdot}$ is the group multiplication of ${\bf G}$ and
$p^{-1}$ is the inverse of $p$. It is easy to see
that $d_{\bf G}$ satisfies the invariance property
$$d_{\bf G}(z\cdot x,z\cdot y)=d_{\bf G}(x,y),  \ \ \forall x, \ y, \ z\in G,\eqno(1.5)$$
and is of homogeneous of degree one, i.e.
$$d_{\bf G}(\delta_\lam(p), \delta_\lam(q))=\lam d_{\bf G}(p,q),\ \forall \lam>0,
\ \forall p, \ q\in\bf G \eqno(1.6)$$
where $\delta_\lam(p)=\lam\xi_1(p)+\sum_{i=2}^r \lam ^i\xi_i(p)$ is the non-isotropic
dilations on $\bf G$. 

Throughout this paper, we fix some notations. For $l\ge 1$,
denote ${\Cal S}(l)$ as the set of $l\times l$ symmetric matrices. 
For $M, N\in {\Cal S}(m)$,
we say $M\ge N$ if $(M-N)\in {\Cal S}(m)$ is a positive semidefinite matrix,
and  let $\hbox{trace}(M)$ denote the trace of $M$ for $M\in{\Cal S}(m)$. 
For $u:{\bf G}\to R$, 
let $\nabla u, \nabla^2 u$ denote the Euclidean gradient, hessian of $u$
respectively, and $\nabla_h u:=(X_1 u, \cdots, X_m u)$, 
$\nabla^2_h u:=({X_iX_j +X_jX_i\over 2}u)_{1\le i, j\le m}$
denote the horizontal gradient, horizontal hessian of $u$
respectively. 
For a given domain $\Om\ci\bf G$, denote $C({\Om})$ 
as the set of continuous functions on $\Om$,
$C^2(\Om)=\{u\in C(\Om): \nabla u, \nabla^2 u\in C(\Om)\}$,
and $\Gamma^2(\Om)=\{u\in C(\Om): \nabla_h u, \nabla^2_h u\in C(\Om)\}$.
A fully nonlinear partial horizontal-differential operator 
${\Cal F}[\cdot]$ on $\Om$ is defined by
$${\Cal F}[\phi](x)=F(\phi(x),\nabla_h\phi(x), \nabla^2_h \phi(x)), \ \forall x\in\Om,
\ \forall \phi\in {\Gamma}^2(\Om), \eqno(1.7)$$
where $F\in C(R\times R^m\times {\Cal S}(m))$. We now give the
definition of subellipticity and nondecreasing property of ${\Cal F}$.
\ss
\nind{\bf Definition 1.1}. The operator ${\Cal F}[\cdot]$ is degenerate subelliptic
if
$$F(r,p, M)\le F(r,p, N), \ \hbox{for all } M, N\in {\Cal S}(m) \hbox{ with } M\le N
\hbox{ and } (r,p)\in R\times R^m. \eqno(1.8)$$
The operator ${\Cal F}[\cdot]$ is uniformly subelliptic if  there exist constants
$\alpha_1, \alpha_2>0$ such that
$$F(r,p,M)-F(r,q,N)\ge \alpha_1\hbox{trace}(M-N)-\alpha_2|p-q| \eqno(1.9)$$
for all $M,N\in{\Cal S}(m)$ with $M\ge N$ and $(r,p,q)\in R\times R^m\times R^m$.
\ss
\nind{\bf Definition 1.2}. The operator ${\Cal F}[\cdot]$ is nonincreasing
if
$$F(r,p,M)\le F(s,p,M),\ \hbox{ for all } r\ge s, \hbox{ and } (p, M)\in R^m\times {\Cal S}(m). 
\eqno(1.10)$$
The operator ${\Cal F}[\cdot]$ is decreasing if there is a constant $\alpha_3>0$ such that
$$F(r,p,M)-F(s,p,M)\le \alpha_3(s-t) \hbox{ for all } r\ge s, \hbox{ and }
(p,M)\in R^m\times{\Cal S}(m). \eqno(1.11)$$

We shall now recall the definition of viscosity solution of fully nonlinear 
degenerate subelliptic equation (1.7),
which was introduced by Crandall-Lions (see [CL] and the user's guides [CIL]) 
for fully nonlinear elliptic equasions.
\ss
\nind{\bf Definition 1.3}. Assume that ${\Cal F}$ is a degenerate subelliptic
operator. $w\in C(\Om)$ is a viscosity subsolution of (1.7) if
for any $(x_0,\phi)\in \Om\times C^2(\Om)$ such that
$$0=\phi(x_0)-w(x_0)\ge \phi(x)-w(x), \forall x\in\Om$$
we have
$$F(\phi(x_0), \nabla_h\phi(x_0), \nabla^2_h \phi(x_0))\ge 0. \eqno(1.12)$$
$w\in C(\Om)$ is a viscosity supersolution of (1.7) if $-w$ is a viscosity
subsolution of (1.7). $w\in C(\Om)$ is a viscosity solution of (1.7) if
it is both a viscosity subsolution and a viscosity supersolution.
\ss
\nind{\bf Remark 1.4}. 
(i) It is well-known that there is an equivalent formulation of viscosity
solution of (1.7) using elliptic jets (see [CIL] or [J1]\S2).
(ii) From the sub-Riemannian point of views, it also
seems natural to define an intrinsic version of viscosity solution of (1.7) by
allowing the test functions $\phi\in \Gamma^2(\Om)$ in Definition 1.3.
However, since $C^2(\Om)\ci \Gamma^2(\Om)$, the intrinsic version of
viscosity solution of (1.7) is stronger than the version given by Definition 1.3.
(iii) This intrinsic version of viscosity solution of (1.7) has been previously
formulated by Bieske [B] (see also Manfredi [M]) in connections with
$\fy$-sublaplacian equations on Heisenberg groups, 
where the notion of subelliptic jets was also introduced.

Now we are ready to state our comparison theorem.
\ss
\nind{\bf Theorem A}. {\it Let $\bf G$ be a Carnot group and $\Om\ci\bf G$ be
a bounded domain. Suppose that $u\in C(\bar\Om)$ is a viscosity subsolution of (1.7)
and $v\in C(\bar\Om)$ is a viscosity supersolution of (1.7). If $F$ satisfies either

(i) ${\Cal F}[\cdot]$ is degenerate subelliptic and decreasing,

\nind or

(ii) ${\Cal F}[\cdot]$ is uniformly subelliptic and nonincreasing,

\nind then
$$\sup_{\Om}(u-v)^+\le \sup_{\del\Om}(u-v)^+. \eqno(1.13)$$}

We would like to remark that the operator $\bar{\Cal F}[\cdot]$ induced by
the degenerate subelliptic operator ${\Cal F}[\cdot]$:
$${\bar{\Cal F}}[w]={\bar F}(x, w, \nabla w, \nabla^2 w):=F(w, \nabla_h w,\nabla^2_h 
w)\eqno(1.14)$$
may not be degenerate elliptic (i.e ${\bar F}$ may not be monotone
in its third variable, see [CIL]), and may be dependent of the spatial variable $x$ in an
essential way so that the uniqueness theorems by Jensen [J1] and Ishii [I]
on viscosity solutions to 2nd order elliptic equations are not 
applicable here. Therefore theorem A not only provides a comparison
principle in the subelliptic setting of the Carnot group but also
makes the comparison principle of [J1] (see also [I] or [CIL])  available
for a  considerably larger class of equations in the Euclidean setting. 

We believe that theorem A shall play an important role in the existence of viscosity
of (1.7) by the Perron's method (see [I]) and plan to study
it in a future article. We would like to mention that Manfredi [M] proved, among other things, 
theorem A for any uniformly elliptic, linear subelliptic opertor 
$F(w)=\sum_{i, j=1}^m a_{ij}(x)X_iX_jw$, with $(a_{ij}(x))\in C({\bf G}, {\Cal S}(m))$
a uniformly elliptic matrix. 

A direct consequence of Theorem A is the uniqueness theorem of
viscosity solutions of (1.7).

\nind{\bf Corollary B}. {\it Under the same assumptions as Theorem A. 
There exists at most one viscosity solution $u\in C(\bar\Om)$ of (1.7).}

The basic point to prove theorem A is that we can always compare between
a {\it classcial} subsolution $u\in C^2(\Om)\cap C(\bar\Om)$
and a {\it classical strict} supersolution $v\in C^2(\Om)\cap C(\bar\Om)$
of (1.7), under the assumption that ${\Cal F}[\cdot]$ is degenerate
subelliptic. 

In order to extend this idea to viscosity sub (or super) 
solutions of (1.7), we first establish, in Lemma 2.1 below, that 
under the same conditions of theorem A, any viscosity
supersolution $v\in C(\bar\Om)$ can be perturbed into a viscosity {\it strict} 
supersolution of (1.7). We would like to point out that, even in the Eucliean setting of $R^n$, 
Lemma 2.1 seems to be new and can be used to simplify the proof of Jensen [J1].
Moreover, it seems necessary in the subelliptic setting, since the counterpart
of [J1] Lemma 3.20 is to estimate ${\hbox{trace}(\nabla^2_h w)^-\over |\nabla_h w|}$ from below for
a semiconvex function $w$ on $\bf G$ and may not be available. 

The second ingredient is to approximate viscosity sub (or super) solutions by
semiconvex (or semiconcave) sub (or super) solutions of (1.7). This idea was first introduced by
Jensen in his very important paper [J1] on uniqueness of Lipschitz continuous
viscosity solutions to 2nd order elliptic equations, and Jensen's original approximation scheme was
further simplified by the sup/inf convolution construction by Jensen-Lions-Souganidis [JLS]
in the Euclidean setting. In the subelliptic setting of the Carnot group,
we succeeded, in an earlier paper [W1] where we proved the uniqueness of viscosity
solution to the subelliptic $\fy$-laplacian equation on any Carnot group $\bf G$, 
in extending the sup/inf convolution construction of [JLS] by employing the smooth 
gauge pseudo-norm function $d_{\bf G}$ and get the desired approximations.  
For the reader's convenience, we review the sup/inf convolution construction
of [W1] in \S3 below. We would like to point out that we have used in a very
crucial way that (1.7) is invariant under  
the group multiplication from left, i.e. if $u\in C({\bf G})$ is
a viscosity solution to (1.7) then $u_a(x)=u(a\cdot x), \ x\in {\bf G},$
is also a viscosity solution to (1.7) for any $a\in{\bf G}$. 
Once we have semiconvex (or semiconcave) sub (or super) solutions to (1.7), 
we can apply both the regularity properties (see Evans-Gariepy [EG]) and Jensen's maximal principle
for semiconvex functions (see [J1]) in our setting. 

The paper is written as follows. In \S2, we show that any viscosity supersolution of (1.7)
given by theorem A can be perturbed into a {\it strict} supersolution of (1.7).
In \S3, we recall the sup/inf convolution construction on a Carnot group $\bf G$,
which was carried out in an earlier paper [W1].
In \S4, we give a proof of theorem A.

\bs
\nind\S2 Viscosity strict supersolutions
\ss
In this section, we show that any viscosity supersolution given by theorem A 
can be perturbed into a viscosity strict supersolution by a suitable small
perturbation. More precisely, we have
\ss
\nind{\bf Lemma 2.1}. {\it Suppose that $F\in C(R\times R^m\times{\Cal S}(m))$
and $v\in C(\Om)$ is a viscosity supersolution to
$${\Cal F}[w]:=F(w,\nabla_h w, \nabla^2_h w)=0, \hbox{ in }\ \Om, \eqno(2.1) $$
under either (i) ${\Cal F}$ is degenerate subelliptic and decreasing or
(ii) ${\Cal F}$ is uniformly subelliptic, nonincreasing.  Then
for any $\delta\in (0,\delta_0)$ 
there are $c_\delta>0$ and $v^\delta\in C(\Om)$ so that
$$v(x)\le v^\delta(x)\le v(x)+\delta \eqno(2.2)$$ 
and 
$v^\delta$ is a viscosity supersolution to
$$F(w,\nabla_h w, \nabla^2_h w)+c_\delta =0, \hbox{ in }\  \Om. \eqno(2.3)$$}
\ss
\nind{\bf Proof}. We first recall that for $x\in\bf G$ if $(x_1,\cdots,x_m)$ denotes the horizontal
component of its coordinate and $((x_{2,1},\cdots, x_{2,m_2}), \cdots, (x_{r,1},\cdots, x_{r,m_r}))$
denotes the vertical component of its coordinate, then the horizontal vector fields $X_l, 1\le l\le m,$
can be expressed as (see [FS])
$$X_l={\del\over\del x_l}+\sum_{i=2}^r\sum_{j=1}^{m_i}a_{ij}(x){\del\over\del x_{i,j}},
\ 1\le l\le m, \eqno(2.4)$$
where $\{a_{ij}\}$ are smooth on $\bf G$ for $2\le i\le r, 1\le j\le m_i$.
For $k>0$, denote $c_1=\inf_{x\in \Om}x_1\in R$ and
define 
$$\alpha_k(x)=1-{1\over k}e^{-k(x_1+1-c_1)}, \ \forall x\in\Om.$$
Then (2.4) implies that, for any $x\in\Om$, we have
$$X_1\alpha_k(x)=e^{-k(x_1+1-c_1)}, \ X_2\alpha_k(x)=\cdots =X_m \alpha(x)=0, \eqno(2.5)$$
$$\eqalignno{X_{ij}\alpha_k(x)&=-ke^{-k(x_1+1-c_1)}, \hbox{ if } i=j=1, &(2.6)\cr
                              &=0, \ \ \ \ \ \ \ \ \ \ \ \ \ \ \ \ \ \hbox{ otherwise}.\cr}$$  

For any $\delta>0$ and $k\ge 2$ to be chosen later, we consider
$v^\delta(x):=v(x)+\delta \alpha_k(x):\Om\to R$. Since $0\le \alpha_k\le 1$, it is
easy to see that $v^\delta$ satisfies (2.2). We want to show that $v^\delta$ is also
a viscosity supersolution to eqn. (2.4). To do this, let $x_0\in \Om$ and 
$\phi\in C^2(\Om)$ touch $v^\delta$ from below at $x=x_0$, i.e.
$$0=(v+\delta \alpha_k-\phi)(x_0)\ge (v+\delta \alpha_k-\phi)(x), \ \forall x\in \Om.$$
This implies that $\phi-\delta \alpha_k\in C^2(\Om)$ touches $v$ from below at $x_0$.
Since $v$ is a viscosity supersolution to eqn.(2.1), we have
$$F(\phi-\delta \alpha_k,\nabla_h \phi-\delta \nabla_h \alpha_k, 
\nabla_h^2\phi-\delta \nabla_h^2 \alpha_k)(x_0)\le 0. \eqno(2.7)$$
Now we need to show that there exists a $k_0>0$ such that for any $k\ge k_0$ (2.3) is true. 
We proceed it as follows.

\nind{\it Case 1}. {\it ${\Cal F}$ is degenerate subelliptic and decreasing }:

It follows from (2.6) that $\nabla^2_h(-\alpha_k)$ is positive semidefinite. Therefore,
the degenerate subellipticity (1.8) and decreasing property (1.11) of ${\Cal F}$ imply 
$$\eqalignno{&F(\phi-\delta \alpha_k, \nabla_h \phi-\delta \nabla_h \alpha_k, 
\nabla_h^2\phi-\delta \nabla_h^2 \alpha_k)(x)\cr
&\ge F(\phi-\delta \alpha_k, \nabla_h \phi-\delta \nabla_h \alpha_k, \nabla_h^2\phi)(x)\cr
&\ge F(\phi, \nabla_h \phi-\delta \nabla_h \alpha_k, \nabla_h^2\phi)(x)+\alpha_3\delta \alpha_k(x)\cr
&\ge F(\phi,\nabla_h\phi,\nabla_h^2\phi)(x)+\alpha_3\delta \alpha_k(x)
-\om_2(\delta |\nabla_h \alpha_k|(x)), &(2.8)\cr}$$
where $\om_2$ is the modular of continuity of $F$ with respect to it second variable.
Since $\alpha_k(x)\ge {1\over 2}$ and $|\nabla_h \alpha_k|(x)\le {1\over 2k}$, (2.8)
implies
$$F(\phi,\nabla_h\phi, \nabla_h^2 \phi)(x_0)
\le \om_2(\delta k^{-1})-{k_0\delta\over 2}\equiv c_\delta<0, \eqno(2.9)$$
if we choose $k$ so large that $\om_2(\delta k^{-1})\le {k_0\delta\over 4}$.
This verifies (2.4) under the condition (i) of Lemma 2.1.

\nind{\it Case 2}. {\it ${\Cal F}$ is uniformly subelliptic, nonincreasing }:

Since $\nabla_h^2(-\alpha_k)$ is positive semidefinite, the uniform ellipticity (1.9)
and nonincreasing property (1.10) of ${\Cal F}$ imply
$$\eqalignno{&F(\phi-\delta \alpha_k,\nabla_h \phi-\delta \nabla_h \alpha_k, 
\nabla_h^2\phi-\delta \nabla_h^2 \alpha_k)(x)\cr
&\ge F(\phi, \nabla_h \phi-\delta \nabla_h \alpha_k, \nabla_h^2\phi)(x)
+\alpha_1\hbox{trace}(\nabla_h^2 (-\alpha_k))(x)\cr
&\ge F(\phi, \nabla_h \phi-\delta \nabla_h \nabla_k, \nabla_h^2\phi)(x)+\alpha_1ke^{-k(x_1+1-c_1)}\cr
&\ge F(\phi,D_h\phi,D_h^2\phi)(x)+\alpha_1ke^{-k(x_1+1-c_1)}-\alpha_2 \delta |\nabla_h\alpha_k|(x)\cr
&\ge F(\phi, D_h\phi, D_h^2\phi)(x)+e^{-k(x_1+1-c_1)}(\alpha_1 k-\alpha_2\delta). &(2.10)\cr}$$
Therefore if we choose $k\ge {2\alpha_2\delta\over k_1}$, then we have
$$F(\phi,\nabla_h\phi,\nabla_h^2\phi)(x_0)\le -{\alpha_2\delta\over 2}e^{-k(x_1+1-c_1)}, \eqno(2.11)$$
this, combined with $\inf_{x\in \Om}e^{-k(x_1+1-c_1)}=2c_2(k)>0$, implies that
(2.4) holds with $c_\delta=c_2(k)\alpha_2\delta>0$. Therefore, the proof of Lemma 2.1 is complete.
\qed

\bs
\nind \S3. The construction of sup/inf convolutions on $\bf G$
\ss
For the convenience of readers, we recall the construction of sup/inf convolution on any
Carnot group $\bf G$, which was carried out earlier by Wang [W1]. The key observation
is that the equation (1.7) is invariant under group multiplication from left on $\bf G$.
We would like to point out that this construction is an extension of that
by Jensen-Lions-Souganidis [JLS] in the Euclidean space.

Let $\Om\ci\bf G$ be a bounded domain and $d_{\bf G}(\cdot,\cdot)$ be
the smooth gauge distance defined by (1.3).  For any $\epsilon>0$, define
$$\Om_\epsilon=\{x\in \Om: \inf_{y\in {\bf G}\setminus\Om} d_{\bf G}(x^{-1},y^{-1})
\ge\epsilon\}.$$
\ss
\nind{\bf Definition 3.1}. For any $u\in C(\bar\Om)$ and $\epsilon>0$, the
sup involution $u_\epsilon$ of $u$ is defined by
$$u^\epsilon(x)=\sup_{y\in\bar\Om}(u(y)-{1\over 2\epsilon}d_{\bf G}(x^{-1},y^{-1})^{2r!}), 
\ \forall x\in\Om. \eqno(3.1)$$
Similarly, the inf involution $v_\epsilon$ of $v\in C(\bar\Om)$ is defined by
$$v_\epsilon(x)=\inf_{y\in\bar\Om}(v(y)+{1\over 2\epsilon}d_{\bf G}(x^{-1},y^{-1})^{2r!}), 
\ \forall x\in\Om. \eqno(3.2)$$

For $p\in \bf G$, let $\|p\|_{E}:=(\sum_{i=1}^r|\xi_i(p)|^2)^{1\over 2}$ 
be the euclidean norm of $p$. We recall
\ss
\nind{\bf Definition 3.2}. A function $u\in C(\bar\Om)$ is called semiconvex, if
there is a constant $C>0$ such that $u(p)+C\|p\|_E^2$ is convex in the Euclidean sense;
and $u$ is called semiconcave if $-u$ is semiconvex. Note that, for $u\in C^2(\Om)$, if
$\nabla^2u(p)+C$ is positive semidefinite for any $p\in\Om$, then $u$ is semiconvex.

Now we have
\ss
\nind{\bf Proposition 3.3}. {\it For $u, v\in C(\bar\Om)$, denote
$R_0=\max\{\|u\|_{L^\fy(\Om)},\|v\|_{L^\fy(\Om)}\}$. Then,
for any $\epsilon>0$,  $u^\epsilon, v_\epsilon\in W^{1,\fy}_{\hbox{cc}}(\Om)$
satisfy

(1) $u^\epsilon$ is semiconvex and $v_\epsilon$ is semiconcave.

(2) $u^\epsilon$ is monotonically nondecreasing w.r.t. $\epsilon$
and converges uniformly to $u$ on $\Om$;
and $v_\epsilon$ is monotonically nonincreasing w.r.t. $\epsilon$ and 
converges uniformly to $v$ on $\Om$.

(3) if $u$ (or $v$ respectively) is a viscosity subsolution (or
supersolution respectively) to a degenerate subelliptic
equation:
$$ F(u, \nabla_h u, \nabla^2_h u) = 0 \ \hbox{ in }\ \Om,  \eqno(3.3)$$
where $F\in C(R\times R^m\times S(m))$.
Then $u^\epsilon$ (or $v_\epsilon$) is a viscosity subsolution
(or supersolution respectively) to eqn. (3.3) in $\Om_{2R_0\epsilon}$.}
\ss
\nind{\bf Proof}. Since the proof of $v_\epsilon$ can be done by the same wasy as
that of $u^\epsilon$, it suffices to consider $u^\epsilon$. 
For $\Om\ci\bf G$ is bounded, the formula (1.3) of $d_{\bf G}$ implies
$$C(\Om,d_{\bf G})\equiv \|\nabla^2_x(d_{\bf G}(x^{-1},y^{-1})^{2r!})
\|_{L^\fy(\Om\times\Om)}<\fy.$$
Therefore, for any $y\in\bar\Om$, the full hessian of
$$\tilde {u^\epsilon}(x,y):=u(y)-{1\over 2\epsilon}
d_{\bf G}(x^{-1},y^{-1})^{2r!}+{C(\Om,d_{\bf G})\over 2\epsilon}\|x\|_{E}^2,
\ \forall x\in\Om,$$
is positive semidefinite so that ${\tilde{u^\epsilon}}$ is convex.
Note that the superum for a family of convex functions is still convex, this
implies that
$$u_\epsilon(x)+{C(\Om,d_{\bf G})\over 2\epsilon}\|x\|_E^2=\sup_{y\in\bar\Om}{\tilde {u^\epsilon}}(x,y),
\ \forall x\in\Om$$
is convex so that $u_\epsilon$ is semiconvex. It is well-known that semiconvex functions are
Lipschitz continuous with respect to the euclidean metric (cf. Evans-Gariepy [EG]).
Therefore $u^\epsilon$ is Lipschitz continuous in $\Om$ with respect to $d_{\bf G}$. This gives (1).

For any $\epsilon_1<\epsilon_2$, it is easy to see that 
$u^{\epsilon_1}(x)\le u^{\epsilon_2}(x)$ 
so that $\{u^\epsilon\}$ is monotonically nondecreasing with
respect to $\epsilon$. Observe that for any $x\in\Om$ there
exists a $x_\epsilon\in\bar\Om$ such that
$$u(x)\le u^\epsilon(x)=u(x_\epsilon)-{1\over 2\epsilon}d_{\bf G}(x^{-1},x_\epsilon^{-1})\le R_0.
\eqno(3.4)$$
This implies
$$u(x_\epsilon)-u^\epsilon(x)={1\over 2\epsilon}d_{\bf G}(x^{-1},x_\epsilon^{-1})
\le u(x_\epsilon)-u(x)=\om_u(\|x_\epsilon-x\|_E),
\ \forall x\in\Om, \eqno(3.5)$$
where $\om_u$ denotes the modular of continuity of $u$.
On the other hand, the monotonicity of $u^\epsilon$ with
respect to $\epsilon$ implies
$$u_{\epsilon\over 2}(x)\ge u(x_\epsilon)-{1\over\epsilon}d(x^{-1}, x_\epsilon^{-1})^{2r!}
=u_\epsilon(x)-{1\over 2\epsilon}d(x^{-1}, x_\epsilon^{-1})^{2r!}$$
so that
$$\lim_{\epsilon\rightarrow 0}{1\over\epsilon}d(x^{-1}, x_\epsilon^{-1})^{2r!}=0,
\ \forall x\in\Om. \eqno(3.6)$$
This implies that $\lim_{\epsilon\rightarrow 0}x_\epsilon=x$ and
$\lim_{\epsilon\rightarrow 0}u^\epsilon(x)=u(x)$ for any $x\in\Om$. 
Moreover, (3.5) implies
$$d_{\bf G}(x^{-1},x_\epsilon^{-1})\le 2\epsilon\om_u(\|x_\epsilon-x\|_E)\le 2R_0\epsilon \eqno(3.7)$$
so that $\|x_\epsilon-x\|_E\le C(R_0\epsilon)^{1\over r}$, where $r$ is the step of $\bf G$.
This, combined with (3.5) again, implies
$$\max_{x\in\Om}|u^\epsilon(x)-u(x)|\le \om_u(C(R_0\epsilon)^{1\over r})\rightarrow 0, 
\ \hbox{ as } \epsilon\rightarrow 0$$
so that $u^\epsilon$ converges to $u$ uniformly. Therefore (2) is proved.

For (3), we first observe that (3.7) implies that for $x^0\in \Om_{2R_0\epsilon}$,
$u^\epsilon(x^0)$ is attained by a $x_\epsilon^0\in\Om$.
Now we let $\phi\in C^2(\Om)$ be such that
$$u^\epsilon(x^0)-\phi(x^0)\ge u^\epsilon(x)-\phi(x), \ \ \forall x\in\Om_{2R_0\epsilon}.$$
Then we have, for any $x, y\in \Om_{2R_0\epsilon}$,
$$u(x^0_\epsilon)-{1\over 2\epsilon}d_{\bf G}((x^0)^{-1},(x^0_\epsilon)^{-1})^{2r!}-\phi(x^0)
\ge u(y)-{1\over 2\epsilon}d_{\bf G}(x^{-1}, y^{-1})^{2r!}-\phi(x). \eqno(3.8)$$
For $y$ near $x^0_\epsilon$, since $x=x^0\cdot (x^0_\epsilon)^{-1}\cdot
y\in \Om_{2R_0\epsilon}$, we can substitue $x$ into (3.8) to get
$$u(x^0_\epsilon)-\phi(x^0)\ge u(y)-\phi(x^0\cdot (x^0_\epsilon)^{-1}\cdot y).$$
Set $\bar{\phi}(y)=\phi(x^0\cdot (x^0_\epsilon)^{-1}\cdot y)$ for $y\in\Om_{2R_0\epsilon}$
close to $y_0$. Then ${\bar\phi}$ touches $u$ from above at $y=x^0_\epsilon$ so that
$u$ being a viscosity subsolution of eqn. (3.3) implies 
$$F(u(x_\epsilon^0), \nabla_h{\bar\phi}(x_\epsilon^0), \nabla^2_h{\bar\phi}(x_\epsilon^0))
\le 0. \eqno(3.9)$$
Note that the left-invariance of $X_i$, we know
$$\nabla_h{\bar\phi}(y)=\nabla_h\phi(x^0\cdot (x^0_\epsilon)^{-1}\cdot y),
\ \  \nabla_h^2{\bar\phi}(y)
=\nabla^2_h \phi(x^0\cdot (x^0_\epsilon)^{-1}\cdot y).$$
Hence we have
$$F(u(x_\epsilon^0), \nabla_h\phi(x_0), \nabla^2_h\phi(x_0))\le 0. \eqno(3.10)$$
Taking $\epsilon$ into zero, (3.10) implies that
$u^\epsilon$ is a viscosity subsolution of eqn.(3.3) on $\Om_{2R_0\epsilon}$.
The proof is complete.           \qed
\bs
\nind \S4. Proof of Theorem A
\ss
This section is devoted to the proof of the comparison Theorem. 
The idea is to prove the comparison property between 
the strict supersolution obtained by Lemma 2.1 and the subsolution
by comparing their sup/inf convolutions. The almost everywhere twice 
differentiablity ([EG]) and Jensen's maximum principle ([J1,2])
for semiconvex functions play very important roles in this aspect.

Through this section, we express the horizontal vector fields $X_i$
by the formula (2.4). 
\ss
\nind{\bf Proof of Theorem A}.

Suppose that (1.13) were fasle.  Then
$$\delta_0=\sup_{\bar\Om}(u-v)^+-\sup_{\del\Om}(u-v)^+>0. \eqno(4.1) $$
Denote $c^+=\sup_{\del\Om}(u-v)^+\ge 0$. Note that $v+c^+$ 
is also a viscosity supersolution to eqn.(?), and (4.1) implies
$$\delta_0=\sup_{\bar\Om}(u-(v+c^+))^+-\sup_{\del\Om}(u-(v+c^+))^+>0.$$ 
Hence we may assume $c^+=0$ (i.e. $u(x)\le v(x)$ for any $x\in\del\Om$)
so that (4.1) implies $\delta_0=\sup_{\Om}(u-v)>0$.

For any $\delta\in (0,{\delta_0\over 4})$, let $v^\delta\in C(\bar\Om)$ be the
strict supersolution of eqn.(1.7) given by Lemma 2.1. In particular,
$v^\delta$ is a viscosity supersoltution to
$$F(w,\nabla_h w, \nabla^2_h w)+c_\delta=0, \ \hbox {in }\ \Om. \eqno(4.2)$$
For any $\epsilon\in (0,\delta)$, we now let $u^\epsilon \ (v^\delta_\epsilon$, respectively$)$
be the sup-convolution (inf-convolution, respectively) of $u \ (v^\delta$ respectively$)$ given
by Proposition 3.3. By considering a smaller domain, we may assume that
$u^\epsilon$ is a viscosity subsolution of eqn.(1.7) and $v^\delta_\epsilon$
is a viscosity supersolution of eqn.(4.2) in $\Om$, and
$$\sup_{\bar\Om}(u^\epsilon-v^\delta_\epsilon)>0\ge\sup_{\del\Om}(u^\epsilon-v^\delta_\epsilon)$$
is achieved at a point $x_0\in\Om$. Since Proposition 3.3 implies
that $u^\epsilon-v^\delta_\epsilon$ is semiconvex, we know (cf. [J2] page 67) that
$$\nabla u^\epsilon(x_0), \nabla v^\delta_\epsilon(x_0) \hbox{ both exist and are equal}, \eqno(4.3)$$
$$\eqalignno{u^\epsilon(x)-u^\epsilon(x_0)-\lan \nabla u^\epsilon(x_0), 
x-x_0\ran_E &=O(\|x-x_0\|_E^2), &(4.4)\cr
v^\delta_\epsilon(x)-v^\delta_\epsilon(x_0)
-\lan \nabla v(x_0), x-x_0\ran_E &=O(\|x-x_0\|_E^2), &(4.5)\cr}$$
where $\lan\cdot, \cdot\ran_E$ denotes the Euclidean inner product on $\bf G$.
Let $R_0=\hbox{dist}_E(x_0,\del\Om)=\inf_{x\in\del\Om}\|x_0-x\|_E>0$ 
be the euclidean distance from
$x_0$ to $\del\Om$ and $R_1>0$ be such that both (4.4) and (4.5) hold with $\|x-x_0\|<R_1$.
Set $R_2=\min\{R_0,R_1\}>0$. 
For simplicity, we will denote $u$, $v$ as $u^\epsilon$, $v^\delta_\epsilon$ respectively
from now on. For any small $\rho>0$, 
define the rescaled maps $u^\rho, v^\rho$ in the euclidean ball $B_{R_2\rho^{-1}}^E$ by
$$\eqalignno{u^\rho(x)&={1\over\rho^2}(u(x_0+\rho x)-u(x_0)-\rho\lan \nabla u(x_0), x\ran_E),\cr
v^\rho(x)&={1\over\rho^2}(v(x_0+\rho x)-v(x_0)-\rho\lan \nabla v(x_0), x\ran_E),\cr}$$
where we have used the Euclidean addition and scalar multiplication.
Then it is easy to see 
$$(u^\rho-v^\rho)(0)>0, \ 
(u^\rho-v^\rho)(0)\ge (u^\rho-v^\rho)(x), \  \forall x\in B_{R_2 \rho^{-1}}^E. \eqno(4.6)$$
It follows from (4.4) and (4.5) that, for any $R>0$, 
there exists an $\rho_0=\rho_0(R)>0$ such that (i)
$\{u^\rho\}_{\{0<\rho\le \rho_0\}}$ are uniformly bounded,
uniformly semiconvex, and uniformly Lipschitz continuous in $B_R^E$;
and (ii) $\{v^\rho\}_{\{0<\rho\le \rho_0\}}$ are uniformly bounded,
uniformly semiconcave, and uniformly Lipschitz continuous in $B_R^E$.
Therefore, by the Cauchy diagonal process, we may assume that
there is $\rho_i\downarrow 0$ such that $u^{\rho_i}\rightarrow u^*$,
$v^{\rho_i}\rightarrow v^*$ locally uniformly in $R^n$, where $n=\hbox{dim}(\bf G)$.
In particular, (i) and (ii) imply that $u^*$ is locally bounded, semiconvex in $R^n$, and
$v^*$ is locally bounded, semiconcave in $R^n$, and
$$(u^*-v^*)(0)>0, \ (u^*-v^*)(0)\ge (u^*-v^*)(x), \ \ \ \forall x\in R^n.$$
Now we have 
\ss
\nind{\bf Claim 4.1}. {\it $u^*$ satisfies, in the sense of viscosity, 
$$F(u(x_0), \nabla_h u(x_0), \{\sum_{k,l=1}^n (A_{kl}^{ij}{\del^2 u^*\over\del x_k\del x_l}+
a_{ik}(x_0){\del a_{jl}\over\del x_k}(x_0){\del u\over\del x_l}(x_0)\}_{1\le i,j\le m})\ge 0, 
\hbox{ in }R^n, \eqno(4.7)$$
and $v^*$ satisfies, in the sense of viscosity,
$$F(v(x_0), \nabla_h v(x_0), \{\sum_{k,l=1}^n (A_{kl}^{ij}{\del^2 v^*\over\del x_k\del x_l}+
a_{ik}(x_0){\del a_{jl}\over\del x_k}(x_0){\del v\over\del x_l}(x_0)\}_{1\le i,j\le m})+c_\delta\le 0,
\hbox{ in }R^n, \eqno(4.8)$$ 
where $A_{kl}^{ij}=a_{ik}(x_0)a_{jl}(x_0)$, for $1\le k, l\le n, 1\le i, j\le m$, and $c_\delta>0$.}

Let's assume Claim 4.1 for the moment and proceed as follows.
Since $u^*-v^*$ is semiconvex and achieves its maximum at $x=0$,
we can apply Jensen's maximum principle for semiconvex functions
(see [J1] [J2]) to conclude that there exists $x_*\in R^n$ such that
$\nabla^2u^*(x_*), \nabla^2v^*(x_*)$ both exist and $\nabla^2(u^*-v^*)(x_*)$ is negative
semidefinite.
Denote $M_1, M_2\in {\Cal S}(m)$ by
$$M_1^{ij}=\sum_{k,l=1}^n \{A_{kl}^{ij}{\del^2 u^*\over\del x_k\del x_l}(x_*)
+a_{ik}(x_0){\del a_{jl}\over\del x_k}(x_0){\del u\over\del x_l}(x_0)\},\  1\le i,j \le m, $$
and
$$M_2^{ij}=\sum_{k,l=1}^n \{A_{kl}^{ij}{\del^2 v^*\over\del x_k\del x_l}(x_*)
+a_{ik}(x_0){\del a_{jl}\over\del x_k}(x_0){\del v\over\del x_l}(x_0)\},\ 1\le i,j \le  m.$$
Since (4.2) implies $\nabla u(x_0)=\nabla v(x_0)$, we have
$$\sum_{1\le i, j\le m}(M_1^{ij}-M_2^{ij})p_ip_j =
\sum_{k,l=1}^n \eta_k\eta_l{\del^2 (u-v)^*\over\del x_k\del x_l}(x_*)\le 0, \ \forall p\in R^m,$$
where $\eta_k=\sum_{i=1}^m p_ia_{ik}(x_0)$, for $1\le k\le n$.
Hence $(M_1-M_2)$ is negative semidefinite. Note also that
$u(x_0)>v(x_0)$, $\nabla_h u(x_0)=\nabla_h v(x_0)$. Therefore
the subellipticity and nonincreasing property of ${\Cal F}$ implies
$$F(u(x_0), \nabla_h u(x_0), M_1)\le F(v(x_0), \nabla_h v(x_0), M_2).\eqno(4.9)$$
This clearly contradicts with (4.7) and (4.8), since (4.7) implies
$$F(u(x_0), \nabla_h u(x_0), M_1)\ge 0, $$
and (4.8) implies
$$F(v(x_0), \nabla_h v(x_0), M_2)\le -c_\delta.$$
Therefore the theorem is proved.

Now we indicate the proof of claim 4.1.
This claim follows from the compactness theorem (cf. [CIL]) among a family
of viscosity sub/supersolutions to 2nd order PDEs. For simplicity, we only indicate
how to prove (4.7).
First we claim that $u^\rho$ satisfies, in the sense of viscosity, in $B_{R_2\rho^{-1}}^E$,
$$\eqalignno{&F(u(x_0+\rho x), \nabla_h u(x_0)+\rho X^\rho(x) u^\rho(x),\cr
&\{\sum_{k,l=1}^n A_{ij,kl}^\rho(x){\del^2 u^\rho\over\del x_k\del x_l}(x)
+B_{ij,l}^\rho(x)({\del u\over\del x_l}(x_0)
+\rho {\del u^\rho\over\del x_l}(x))\}_{1\le i,j \le m})=0, &(4.10)\cr}$$
where $A_{ij,kl}^\rho(x) =a^\rho_{ik}a^\rho_{jl}(x)$,
$B_{ij,l}^\rho(x)=\sum_{k=1}^n a^\rho_{ik}({\del a_{jl}\over\del x_k})^\rho$,
$X^\rho(x)=(X_1^\rho(x),\cdots, X_m^\rho(x))$, $X_i^\rho(x)=X_i(x_0+\rho x)$,
$a_{ik}^\rho(x)=a_{ik}(x_0+\rho x)$, and
$({\del a_{jl}\over\del x_k})^\rho(x)={\del a_{jl}\over\del x_k}(x_0+\rho x)$.

To see (4.10), let $({\bar x},\phi)\in B_{R_2\rho^{-1}}^E\times C^2(B_{R_2\rho^{-1}}^E)$ be such that
$$0=u^\rho({\bar x})-\phi({\bar x})\ge u^\rho(x)-\phi(x), \ \forall x\in B_{R_2\rho^{-1}}^E.$$
It is straightforward to see 
$$\phi_\rho(x)\equiv u(x_0)+\lan \nabla u(x_0), x-x_0\ran_E +\rho^2\phi({x-x_0\over\rho}),
\forall x\in B_{R_2}^E(x_0) $$ satisfies
$$0=u(x_0+\rho {\bar x})-\phi_\rho(x_0+\rho {\bar x})
\ge u(x)-\phi_\rho(x), \ \forall x\in B_{R_2}^E(x_0).$$
This, combined with the fact that $u$ is a viscosity subsolution to eqn.(1.?),
implies
$$F(u(x_0+\rho {\bar x}), \nabla_h\phi_\rho(x_0+\rho{\bar x}), 
\nabla^2_h\phi_\rho(x_0+\rho {\bar x}))\ge 0. \eqno(4.11)$$
Direct calculations yield
$${\del \phi_\rho\over\del x_k}(x_0+\rho{\bar x})={\del u\over\del x_k}(x_0)
+\rho {\del\phi\over\del x_k}({\bar x}), \ \forall 1\le k\le n,$$
$${\del^2\phi_\rho\over\del x_k \del x_l}(x_0+\rho{\bar x})
={\del^2\phi\over\del x_k \del x_l}({\bar x}), \ \forall 1\le k, l\le n.$$
Substituting these into (4.11), we obtain (4.10).

It is clear that, by taking $\rho\rightarrow 0$, (4.10) implies (4.7).
This proves claim 4.1.                 \qed

\bs
\bs
\cl{\bf REFERENCES}
 
\ms
\nind{[A]} G. Aronsson, {\it Extension of functions satisfying Lipschitz conditions}.
Ark. Mat. 6 (1967), 551-561.

\nind{[B]} T. Bieske, {\it On $\infty$-harmonic functions on the Heisenberg group}. 
Comm. Partial Differential Equations 27 (2002), no. 3-4, 727--761.

\nind{[BC]} T. Bieske, L. Capogna, {\it The Aronsson-Euler equation for absolute
minimizing Lipschitz extensions with respect to Carnot-Carath\'edory metrics}. Preprint (2002).

\nind{[BR]} Z. Balogh, M. Rickly, {\it Regularity of convex functions on Heisenberg groups}.
Preprint.

\nind{[CC]} X. Cabr\'e, L. Caffarelli, {Fully nonlinear elliptic equations}. AMS colloquium publications
43, AMS, Providence, RI, 1995.

\nind{[CIL]} M. Crandall, H. Ishii, P. L. Lions, {\it User's guide to viscosity solutions
of second order partial differential equations}.
Bull. Amer. Math. Soc. (N.S.) 27 (1992), no. 1, 1--67.

\nind{[CL]} M. Crandall, P. L. Lions, {\it Viscosity solutions of Hamilton-Jacobi equations}. 
Trans. Amer. Math. Soc. 277 (1983), no. 1, 1--42. 

\nind{[DGN]} D. Danielli, N. Garofalo, D. Nhieu, {\it Notations of convexity in Carnot groups}.
Comm. Anal. Geom., to appear.

\nind{[E1]} L. Evans, {\it A convergence theorem for solutions of nonlinear second order
elliptic equations}. Indiana Univ. Math. J. 27 (1978), 875-887.

\nind{[E2]} L. Evans, {\it On solving certain nonlinear partial differential equations by
accretive operator methods}. Israel J. Math. 36 (1980), 225-247.

\nind{[EG]} L. Evans, R. Gariepy, {Measure theory and fine properties of functions}. CRC Press, 1992.

\nind{[FS]} G. Folland, E. Stein, Hardy spaces on homogeneous groups. Mathematical Notes, 28. 
Princeton University Press, Princeton, N.J., 1982.

\nind{[G]} C. Gutierrez, {The Monge-Amp\'ere equation}. Progress in nonlinear differential equations
and theire applications, 44, Birkhauser, Boston, MA, 2001.

\nind{[GM]} C. Gutierrez, A. Montanari, {\it Maximum and comparison principles for convex functions
on Heisenberg group}. Preprint.

\nind{[GT]} N. Garofalo, F. Tournier, {\it Monotonicity and estimates of the supremum for
Monge-Amp\'ere measures in the Heisenberg group}. Preprint.

\nind{[I]} H. Ishii, {\it On existence and uniqueness of viscosity 
solutions of fully nonlinear second-order elliptic PDEs}. Comm. Pure Appl. Math.
42 (1989) 14-45.

\nind{[J1]} R. Jensen,  {\it The maximum principle for viscosity solutions of fully nonlinear
second order partial differential equations}. Arch. Rational Mech. Anal. 101 (1988), no. 1, 1--27. 

\nind{[J2]} R. Jensen, {\it Uniqueness of Lipschitz extensions: minimizing
the sup norm of the gradient}. Arch. Rational Mech. Anal. 123 (1993), no. 1, 51--74.   

\nind{[JLS]} R. Jensen,  P. L. Lions, P. Souganidis, {\it A uniqueness result for viscosity
solutions of second order fully nonlinear partial differential equations}. 
Proc. Amer. Math. Soc. 102 (1988), no. 4, 975--978. 

\nind{[JM]} P. Juutinen, J. Manfredi. In preparation.

\nind{[LMS]} G. Lu, J. Manfredi, B. Stroffolini, {\it Convex functions on Heisenberg group}.
Calc. Var., to appear. 

\nind{[M]} J. Manfredi, {\it Fully nonlinear subelliptic equations}. In preparation.

\nind{[W1]} C. Y. Wang, {\it The Aronsson equation for absolute minimizers of $L^\fy$-functionals
associated with vector fields satisfying H\"ormander's condition}. Preprint (2003), avaiable
at http//arXiv: math.AP/0307198.

\nind{[W2]} C. Y. Wang, {\it Viscosity convex functions on Carnot groups}. Preprint (2003).

\end